\documentclass[12pt]{article}
\usepackage{latexsym}
\usepackage{amsmath,amsfonts,amssymb}
\usepackage{epsfig}

\begin{document}
%\parskip\baselineskip
%\makeatletter
%\if\c@secnumdepth>0
\newtheorem{definition}{Definition}[section]
%\else
%\newtheorem{definition}{Definition}
%\fi
%\makeatother
%\newtheorem{definition}{Definition}[subsection]
\newtheorem{theorem}[definition]{Theorem}
\newtheorem{lemma}[definition]{Lemma}
\newtheorem{proposition}[definition]{Proposition}
\newtheorem{examples}[definition]{Examples}
\newtheorem{corollary}[definition]{Corollary}
\def\square{\Box}
%\newremark{example}[definition]{Example}
%\newremark{examples}[definition]{Examples}
\newtheorem{remark}[definition]{Remark}
\newtheorem{remarks}[definition]{Remarks}
\newtheorem{exercise}[definition]{Exercise}
\newtheorem{example}[definition]{Example}
\newtheorem{observation}[definition]{Observation}
\newtheorem{observations}[definition]{Observations}
\newtheorem{algorithm}[definition]{Algorithm}
\newtheorem{criterion}[definition]{Criterion}
\newtheorem{algcrit}[definition]{Algorithm and criterion}

\pagestyle{plain}

%\newremark{assumption}[definition]{Assumption}
\newenvironment{prf}[1]{\trivlist
\item[\hskip \labelsep{\it
#1.\hspace*{.3em}}]}{~\hspace{\fill}~$\square$\endtrivlist}

\title{Stokes matrices for the quantum differential equations of some Fano varieties} 
\author{John Alexander Cruz Morales\footnote{e-mail: alekosandro@gmail.com, jacruzm@impa.br}\\
Instituto Nacional de Matem\'{a}tica Pura e Aplicada, IMPA. \\
Estrada Dona Castorina 110. Rio de Janeiro, 22460-320, Brasil.\\
Marius van der Put\footnote{e-mail: m.van.der.put@rug.nl}\\
Department of Mathematics.\\ 
University of Groningen.\\
9700 AK Groningen, P.O. Box 407, The Netherlands.}
%\author{John Alexander Cruz Morales and Marius van der Put}    
%\address{Tokyo Metropolitan University, Faculty of Science and Engineering, Department of Mathematics and Information Sciences, Minami-Ohsawa 1-1, Hachioji, Tokyo 192-037, Japan.}
%\email{cruzmorales-johnalexander@ed.tmu.ac.jp, alekosandro@gmail.com} 
%\address{University of Groningen, Department of Mathematics P.O. Box 407, 9700 AK Groningen, The Netherlands}
%\email{m.van.der.put@rug.nl}
\date{}

\maketitle

\begin{abstract}  
The classical Stokes matrices for the quantum differential equation of $\mathbb{P}^n$ are
computed using multisummation and the `monodromy identity'. Thus, we recover the results of
D.~Guzzetti that confirm Dubrovin's conjecture for projective spaces. The same method
yields explicit formulas for the Stokes matrices of the quantum differential equations of smooth Fano hypersurfaces in $\mathbb{P}^n$ and for weighted projective spaces.\end{abstract}

%\footnote{MSC2000: 34M40, 53D45. keywords: Stokes matrices, quantum cohomology.} %\end{abstract}

\section{Introduction}
For a Fano variety $X$ one can define a Frobenius structure for
its cohomology and the latter induces a linear differential equation 
(or connection in one or more variables) which is called the quantum differential equation of $X$. This equation reflects geometric properties of $X$ and for many 
varieties $X$ the quantum differential 
equation is explicitly known, see \cite{Guest}, \cite{GuestSakai}.  For the cases that we consider, the quantum differential equation is
an ordinary linear differential equation in a complex variable $z$ and has two singular points $z=0$ and
$z=\infty$. The point $z=0$ is regular singular and the point $z=\infty$ is irregular singular. At
$z=\infty$ the difference between formal (symbolic) solutions and actual solutions in sectors is measured by 
{\it Stokes data}. The contribution of this paper to the theory of quantum differential equations is an explicit computation of the Stokes data by means of the formalism of multisummation. This formalism is the work of many experts (see \cite{vdP-Si}, \S 7.1)  and in \S \ref{s2} we will explain how it can be used to compute  the Stokes data in a purely algebraic way. We note that for a general irregular singularity there are only analytic methods for the determination of the Stokes data. Thus quantum differential equations are rather special.\\

In the remaining part of this introduction we sketch, for the convenience of the reader (with many black boxes and without any originality, compare \cite{Guest,Guz,Ueda1,Ueda2}), some of the theory of quantum cohomology. The relation with the above Stokes data and our results  concerning these are presented.\\

Let $X$ be a (smooth) complex projective Fano variety.  Put
$H^*(X,\mathbb{C})=\oplus _{d\geq 0}H^{2d}(X,\mathbb{C})$.
Let $b_1,\dots ,b_r$ be a basis of $H^2(X,\mathbb{C})$. For 
$t=\sum t_ib_i$, one defines a deformation $\circ_t$
of the usual cup product $\circ$ on $H^*(X,\mathbb{C})$. This deformation
is called the {\it  small quantum product}.  One writes formally
$q_i=e^{t_i}$ and $\partial_i=q_i\frac{\partial}{\partial q_i}$. Further,
$\hbar$ will denote a complex parameter. One defines a connection $\nabla$, called the {\it Dubrovin-Givental connection}, on the trivial vector bundle 
$H^2(X,\mathbb{C}) \times H^{\ast}(X,\mathbb{C}) \rightarrow H^2(X,\mathbb{C})$  by the formula  $\nabla_{\partial_i} = \partial_i - \frac{1}{\hbar} b_i \circ_t$ for $i=1,\dots ,r$. The \emph{quantum differential equations} are the equations $\hbar \partial_i \Psi = b_i \circ_t \Psi$ for $i=1,\dots ,r$ and for functions $\Psi: H^2(X,\mathbb{C}) \rightarrow H^{\ast}(X,\mathbb{C})$. \\

Above, we have supposed  $t \in H^2(X,\mathbb{C})$.  However, it is important to  consider also $t \in H^{\ast}(X,\mathbb{C})$. In that case the deformation of the cup product is called the \emph{big quantum product}. For the corresponding `big quantum cohomology and  connection' we refer to \cite{Dub96, Dub99}. \\

In the sequel we restrict ourselves to the small quantum product and
to the case $r=1$, i.e., the case where the quantum differential equation is an ordinary linear differential equation. For a detailed discussion we refer the reader to \cite{Guest} and references therein. \\

A `good Fano variety' $X$ is a Fano variety such that
$D^bcoh(X)$, the derived category of the coherent sheaves on $X$, is generated as triangulated category by an exceptional collection $(\mathcal{E}_i)_{i=1}^N$. An object $\mathcal{E}$ is exceptional if $\rm{Ext}^i(\mathcal{E},\mathcal{E})$ equals $\mathbb{C}$ for $i=0$ and equals $0$ for $i>0$. Further $(\mathcal{E}_i)_{i=1}^N$ is an exceptional  collection if each 
$\mathcal{E}_i$ is exceptional and  ${\rm Ext}^k(\mathcal{E}_i,\mathcal{E}_j)=0$ for any $i>j$ and any $k$. In this situation the {\it Gram matrix} $G$ of $X$ is defined by
$G_{i,j}=\sum _k(-1)^k \dim {\rm Ext}^k(\mathcal{E}_i,\mathcal{E}_j)$. \\

One of {\it conjectures of Dubrovin} (see \cite{Dub98}) states that the Gram matrix of $X$ coincides with the {\it Stokes matrix} of the quantum differential equation of $X$ (up to a certain equivalence which we will make more explicit). For the complex projective space $\mathbb{P}^{n-1}$, the ordered set of line bundles $O,O(1),\dots ,O(n-1)$ is an exceptional collection
and the Gram matrix $G=(G_{i,j})$ is given by $G_{i,j}={n-1+j-i\choose j-i }$ for $i\leq j$ and $G_{i,j}=0$ otherwise. The inverse $(a_{i,j})$ of $G$, which is equivalent to $G$, has the data  $a_{i,j}=(-1)^{j-i}{n\choose j-i}$ for $i\leq j$ and $a_{i,j}=0$ otherwise. \\

Now we will explain the relation between `our' Stokes data and the Stokes matrix considered in quantum cohomology by Dubrovin et al.. The latter we will call `quantum Stokes matrices' and denote by $St_{qc}$. The irregular singularity of the quantum differential equation at $z=\infty$ has Poincar\'e rank 1. This implies that a given formal (or symbolic) fundamental matrix can be lifted to an actual analytic fundamental matrix on a sector at $z=\infty$ of opening slightly larger than $\pi$. Moreover these liftings are unique.   Let $\Phi_{right}$ and  $\Phi_{left}$ denote two of these lifts, then $St_{qc}$ is defined by $\Phi_{right}=\Phi_{left}St_{qc}$. \\

The multisummation theory produces for every singular direction $d$ of the 
differential equation a Stokes matrix, denoted by $St_d$. This expresses the relation
between multisummation of the formal fundamental matrix left and right of the singular direction $d$. One concludes that $St_{qc}$ equals the ordered product $\prod _dSt_d$  taken over the singular directions $d$ in an interval of lenght $\pi$ (in fact $d\in [0,1/2)$ in our notation).  It so happens that each $St_d$ has only one interesting entry. The collection of these entries  will be called the {\it Stokes data}. We note that `our' Stokes data are closely related to what are called `Stokes factors' in \cite{Guz}.  \\

For the complex projective space $\mathbb{P}^{n-1}$ the conjecture of Dubrovin has
been proved by D.~Guzzetti \cite{Guz}. The matrix $St_{qc}$ 
(the product $\prod _dSt_d$) is a unipotent matrix and is, a priori, rather complicated with respect to the given basis (see \S 6 of \cite{Guz}).
This basis is changed (this is the equivalence mentioned before)  by a permutation, by putting signs and the action of a braid group. The quantum differential equation lives in a family
(in fact induced by the big quantum product),
parametrized by $\mathbb{C}^n\setminus$ the diagonals, of similar
equations where the singular directions at $z=\infty$ vary. The braid group action is derived from loops in this family.  Guzzetti showed that $St_{qc}$ has, w.r.t. a new basis and up to signs, the form $(a_{i,j})$ which proves the Dubrovin's conjecture for $\mathbb{P}^{n-1}$.\\

\noindent
{\it Our results, Theorem \eqref{teorema31}, for the Stokes data $\{x_{\ell ,k}\}_{0\leq k,\ell <n;\ k\neq \ell}$ of $\mathbb{P}^{n-1}$ are}:\\ 
 {\it For odd $n$} and $0\leq \ell <k$ one has $x_{\ell ,k}=-(-1)^{k-\ell}{n\choose k-\ell }$ and $x_{\ell ,k}=-x_{k,\ell}$.\\ 
{\it For even $n$} and $0\leq \ell <k$ one has
$x_{\ell ,k}=-(-1)^{k-\ell}{n\choose k-\ell }$ if $k-\ell \leq \frac{n}{2} $ and 
$x_{\ell ,k}=(-1)^{k-\ell}{n\choose k-\ell }$ if
 $k-\ell > \frac{n}{2} $. \\ 
{\it For even $n$} and $0\leq k<\ell$ one has $x_{\ell ,k}=(-1)^{\ell -k}{n\choose \ell -k }$.\\ 

Theorem \eqref{teorema31} proves again Dubrovin's conjecture for $\mathbb{P}^{n-1}$ and we observe  that the above matrix $(a_{i,j})$, equivalent to $St_{qc}$, can rather simply be expressed into the Stokes data $\{x_{\ell ,k}\}$. 
The Stokes data can be read off from the {\it monodromy identity}
which compares the topological monodromy at $z=0$ with
the Stokes matrices $St_d$ and the formal monodromy at $z=\infty$.
The same method leads to the further results:  computations of the Stokes data for weighted projective spaces (Remark \eqref{rm32} and Proposition \eqref{p33}) and for Fano hypersurfaces (Theorem \eqref{teorema41}). \\

Recent papers on the computation of quantum Stokes matrices
are  \cite{Tan} and \cite{Ueda1,Ueda2}.  The first one proposes 
another proof of Dubrovin's conjecture for $\mathbb{P}^n$. 
In the other two papers quantum Stokes matrices are computed for Grassmannians (based on the results for $\mathbb{P}^n$)
and for cubic surfaces.\\

After completing the calculations of this paper we became aware that a related discussion (from a physical point of view) to our work is presented in \cite{Zas}, for the case of projective spaces. However, the argument in loc.cit. concerns the computation of the Stokes matrices for the so-called $tt^{*}$-equations (see \cite{CV}). The question whether these equations are related to the equations for the quantum cohomology and, in particular, whether their Stokes matrices coincide,  is discussed in \cite{Guestetal}. 

\bigskip

The paper is organized as follows. In section \ref{s2} we give a brief presentation of the theory of Stokes matrices emphasizing the relevant facts for our computation. In section \ref{s3} we present the explicit computation for the case of (weighted) projective spaces and in section \ref{s4} we extend that computation to the case of smooth Fano hypersurfaces. In the sequel $q$ will be replaced by $z$ and the parameter $\hbar$ is taken to be 1. We will often write $\delta$ for $z\frac{d}{dz}$. The quantum differential equation in operator form
for $\mathbb{P}^{n-1}$ then obtains the simple form $\delta ^n-z$.

\section{Stokes matrices and  the monodromy identity} \label{s2}

A linear differential operator of order $n$, analytic in the neighbourhood of $z=\infty$, has a scalar form
 $(z\frac{d}{dz})^n+a_{n-1}(z\frac{d}{dz})^{n-1}+\cdots +a_1z\frac{d}{dz}+a_0$ with all $a_j$ in the field $\mathbb{C}(\{z^{-1}\})$ of the convergent Laurent series in $z^{-1}$. The scalar operator can be transformed into a matrix differential operator
$z\frac{d}{dz}+A$ where the entries of the matrix $A$ are in $\mathbb{C}(\{z^{-1}\})$.\\

 As a {\it differential module} over  $\mathbb{C}(\{z^{-1}\})$, the scalar equation above translates into a vector space $M$ of dimension $n$ over this field, equipped with a $\mathbb{C}$-linear operator $\delta_M$ satisfying
 $\delta _M(fm)=z\frac{d}{dz}(f)\cdot m+f\delta _M (m)$, for 
 $f\in  \mathbb{C}(\{z^{-1}\}),\ m\in M$.  Note that for a suitable
basis of $M$, the matrix $A$ above is the matrix  of $\delta_M$ with respect to this basis. \\
 
 The {\it formal classification} of $M$ is the classification of the differential module 
 $\mathbb{C}((z^{-1}))\otimes M$ over the field $\mathbb{C}((z^{-1}))$ of the formal Laurent series in $z^{-1}$. In general, a root $z^{1/m}$ of $z$ for certain $m\geq 1$ is needed for the
formulation of the classification that we describe now. 

There are distinct elements $q_1,\dots ,q_s\in  z^{\frac{1}{m}}\mathbb{C}[z^{\frac{1}{m}}]$, called the {\it generalized eigenvalues}  of $M$ such that $\mathbb{C}((z^{-1/m}))\otimes M$ is a direct sum of (differential) submodules $N_1,\dots ,N_r$ over $\mathbb{C}((z^{-1/m}))$. The differential module $N_j$ has a basis  such that the operator $\delta_{N_j}$ has the form $q_j\cdot id +\ell_j$, where $\ell_j$ has entries in $\mathbb{C}$. The $q_j$ and the decomposition 
$\mathbb{C}((z^{-1/m}))\otimes M=N_1\oplus \cdots \oplus N_r$ are unique.  The $\ell_j$ are not unique.\\

One defines {\it symbols} $z^\lambda$ for every $\lambda \in \mathbb{C}$, $\log z$ and 
$e(q)$ for every $q\in \cup _{n\geq 1}z^{1/n}\mathbb{C}[z^{1/n}]$, by the rules
$z^{\lambda_1+\lambda _2}=z^{\lambda _1}z^{\lambda _2}$, $z^0=1,\ z^1=z$, $e(q_1+q_2)=e(q_1)e(q_2)$, $e(0)=1$ and $\delta(z^\lambda )=\lambda z^\lambda$,
$\delta (\log z)=1$, $\delta (e(q))=q\cdot e(q)$.  On a sector at $z=\infty$ these symbols have
an obvious interpretation (e.g., the interpretation of $e(q)$ is
$e^{\int q\frac{dz}{z}}$), but not on a full neighbourhood of $z=\infty$.

 Let $\gamma$
denote the automorphism of $\cup_{n\geq 1}\mathbb{C}((z^{-1/n}))$ defined by 
$\gamma z^\lambda = e^{2\pi i \lambda} z^\lambda $ for all $\lambda \in \mathbb{Q}$.
The natural action of $\gamma$ on the symbols is given by the formulas 
$\gamma z^\lambda =e^{2\pi i \lambda}z^\lambda$ for all $\lambda \in \mathbb{C}$,
$\gamma \log z=2\pi i+\log z$, $\gamma e(q)=e(\gamma q)$. \\

{\it The symbolic solution space}.
Let $U$ be the $\mathbb{C}((z^{-1}))$-algebra generated by these symbols. Then $U$ is a
universal Picard--Vessiot ring for the differential field $\mathbb{C}((z^{-1}))$, wich means
that for every differential module $M$ over   $\mathbb{C}((z^{-1}))$, the $\mathbb{C}$-vector
space $V:=\ker (\delta ,U\otimes M)$ has the property that the obvious map
$U\otimes _\mathbb{C}V\rightarrow U\otimes M$ is an isomorphism. Moreover, $U$ is minimal
with this property and $U$ has only trivial differential ideals.  The space $V$ is called {\it the 
symbolic solution space of $M$}. Let $b_1,\dots ,b_d$ be a basis of $M$ over $\mathbb{C}((z^{-1}))$. The elements of $V$ are sums $\sum _{j=1}^d \alpha _jb_j$ where the $\alpha _j\in U$
are (by definition) expressions using formal power series, and the symbols $z^\lambda, \log z,
e(q)$.

The decomposition $U=\oplus _{q} U_q$ with $U_q:= e(q)\mathbb{C}((z^{-1}))[\{z^\lambda \}, \log z]$ 
induces a decomposition $V=\oplus _qV_q$ with $V_q=\ker (\delta ,U_q\otimes M)$. Further
$\gamma$ acts as a $\mathbb{C}$-linear automorphism on $V$ and has the property 
$\gamma (V_q)=V_{\gamma q}$. The action of $\gamma$ on $V$ is called {\it the formal
monodromy}.

Thus we have associated to $M$ a tuple $(V, \{V_q\},\gamma )$
of a finite dimensional $\mathbb{C}$-vector space $V$, a subspace $V_q$ for every
 $q$ in the set of generalized eigenvalues $\cup _{n\geq 1}z^{1/n}\mathbb{C}[z^{1/n}]$,
an element $\gamma \in {\rm GL}(V)$, such that $V=\oplus V_q$ and $\gamma V_q=V_{\gamma q}$ for every $q$. This construction yields in fact an equivalence of Tannakian categories (see \cite{vdP-Si} for more details).\\

{\it Singular directions and multisummation}.
For a pair of distinct eigenvalues $(q ,\tilde{q})$, one considers the operator 
 $z\frac{d}{dz}-(q-\tilde{q})= z\frac{d}{dz}- (cz^\lambda +\cdots)$ with $\lambda >0, \ c\neq 0$ and the
 dots are terms $*z^\mu$ with $0<\mu <\lambda$. The solution of the equation is
 $y:=e^{\frac{1}{\lambda}cz^\lambda +\cdots}$. Let  $d\in \mathbb{R}$ stand for the direction
 $e^{2\pi id}$ at $z=\infty$. Then a real number $d$ is called a {\it singular direction for the pair $(q,\tilde{q})$} if and only if $\frac{c}{\lambda}e^{2\pi i\lambda d}$ is real and negative. 

Let $M$ be a differential module over $\mathbb{C}(\{z^{-1}\})$. {\it Multisummation in a direction $d$} is a $\mathbb{C}$-linear bijection $m_d$ from the symbolic solution space $V$ of $M$ to the space of the actual solutions of $M$ in a sector around the direction $d$. The map $m_d$ exists (and is unique) if $d$ is not a singular direction for any pair $(q,\tilde{q})$ of
eigenvalues of $M$.\\

{\it The Stokes maps}. Let a differential module $M$ over $\mathbb{C}(\{z^{-1}\})$ be given
and let $(V,\{V_q\},\gamma)$ be the tuple corresponding to $\mathbb{C}((z^{-1}))\otimes M$. 
Let $d$ be a direction. Then the Stokes map $St_d$ for this direction has the form
$St_d=1+\sum M_{d,q,\tilde{q}}$, where the sum is taken over all pairs $(q,\tilde{q})$ such 
that $V_q,V_{\tilde{q}}\neq 0$ (i.e., $q$ and $\tilde{q}$ are eigenvalues for $M$), $d$ is a singular
direction for $(q,\tilde{q})$ and
 $M_{d,q,\tilde{q}}:V\overset{projection}{\rightarrow}V_q\overset{linear}{\rightarrow}V_{\tilde{q}}
 \overset{inclusion}{\rightarrow}V$. This Stokes map is obtained by comparing the
 multisummation maps $m_{d-\epsilon}, m_{d+\epsilon}$ (with small enough $\epsilon >0$) from $V$ to actual solutions of the differential equation in a sector around the direction $d$. 
Further $\gamma ^{-1}St_d\gamma =St_{d+1}$.
We note that a direction $d$ can be singular for more than one pair $(q,\tilde{q})$. 

For a given
differential module $M$ over $\mathbb{C}(\{z^{-1}\})$, there is an algorithm computing the tuple $(V,\{V_q\},\gamma)$. The entries of the Stokes maps can be expressed as certain involved
integrals and, in general, these cannot be made explicit.\\

Now we have associated to a differential module $M$ over $\mathbb{C}(\{z^{-1}\})$ 
a tuple $(V,\{V_q\},\gamma ,\{St_d\})$ with the properties stated above. This yields an equivalence
between the Tannakian categories of the differential modules over $\mathbb{C}(\{z^{-1}\})$
and the category of these tuples (see Theorem 9.11 in \cite{vdP-Si}).\\

{\it  A change of variables}.  The inclusion $K:=\mathbb{C}(\{z^{-1}\})\rightarrow
K_n:=\mathbb{C}(\{u^{-1}\})$ with $z=u^n$ and $n>1$ induces a functor which associates
to a differential module $M$ over $K$ the differential module $K_n\otimes M$ over $K_n$.
The corresponding morphism between tuples, maps a tuple $(V,\{V_q\},\gamma ,\{St_d\})$ to  
a tuple $(V,\{V_{\tilde{q}}\},\tilde{\gamma},\{\tilde{St}_d\})$. It can be verified that 
  $V_{\tilde{q}}=V_q$ for $\tilde{q}(u)=q(u^n)$, $\tilde{\gamma}=\gamma ^n$ and 
$\tilde{St}_{d}=St_{nd}$.   Using this one can compare  the 
singularities of, for instance, $(z\frac{d}{dz})^n-z$ and $(u\frac{d}{du})^n-n^nu^n$ where $z=u^n$.

\bigskip

{\it The monodromy identity}.
Let the differential module $M$ over $\mathbb{C}(\{z^{-1}\})$ correspond to the tuple
 $(V,\{V_q\},\gamma ,\{St_d\})$. Let $W$ be a solution space at a certain point $p$ close to $z=\infty$. One makes a loop around $z=\infty$ and analytic continuation along this loop yields the topological monodromy  $mon_\infty \in {\rm GL}(W)$. After some  identification of $W$ with $V$ one obtains the {\it  monodromy identity} (see Proposition 8.12 in \cite{vdP-Si}):
\[ mon_\infty \ \ \ \ \ \mbox{ is conjugated to }\ \ \ \ \ \ \ \ \gamma  \prod _{d\in [0,1),\ d \mbox{ singular }} St_d ,\]
where the order of the maps $St_d$ in the product is  counter clockwise. \\
 
\section{The Stokes matrices for  $\delta^n-z$.} \label{s3}
We summarize the results for this quantum differential operator of $\mathbb{P}^{n-1}$ (normalized by puting $\hbar = 1$). 
The irregular singular point $z=\infty$ has (generalized) eigenvalues $q_j=e^{2\pi ij/n}z^{1/n},\ j=0,\dots , n-1$.

 The symbolic solution space $V$ at $z=\infty$ has a  basis $e_0,\dots ,e_{n-1}$, uniquely determined (up to simultaneous multiplication by a constant) by normalizing the matrix of $\gamma$. Let $E_{k,\ell}\in {\rm End}(V)$ denote the map defined by $E_{k,\ell}e_\ell=e_k$ and $E_{k,\ell}e_j=0$ for $j\neq \ell$. 
For a direction $d$, the Stokes matrix $St_d\in {\rm GL}(V)$ has the form 
$St_d=1+\sum x_{\ell,k}E_{\ell,k}$, where the sum is taken over the pairs $(k,\ell )$ such that the direction $d$ is singular for $q_k-q_\ell$.  For $k\neq \ell$ the pair $(q_k,q_\ell)$ has in the interval $[0,n)$ precisely one singular direction and produces the constant $x_{\ell,k}$. 
Of the $n(n-1)$ singular directions in $[0,n)$ (counted with
multipicity) there are $n-1$ in the interval $[0,1)$. The $x_{\ell ,k}$
corresponding to the singular directions in $[0,1)$ are
computed, using the monodromy identity.   The other $x_{\ell ,k}$ are obtained by the formula $\gamma^{-1}St_d\gamma =St_{d+1}$.

The `{\it Stokes data}' for the equation is by definition $\{x_{\ell,k}\}_{k\neq \ell}$. We note that $x_{\ell,k}=x_{\ell',k'}$ if
$\ell \equiv \ell',\ k\equiv k' \mod n$.
%The indices
%$\ell, k$ are regarded as elements of $\mathbb{Z}/n\mathbb{Z}$. 
The result of this section is:
\begin{theorem} \label{teorema31}  The monodromy identity yields the following formulas:\\
{\bf For $n$ odd}\\
$x_{l,k}=-(-1)^{k-l}{n\choose k-l }\mbox{ for }n>k>l\geq 0 \mbox{ and } k+l=[\frac{n}{2}]\mbox{ or }=[\frac{n}{2}]-1,$\\
$x_{l,k}=(-1)^{l-k}{n\choose l-k }\mbox{ for }n> l>k\geq 0 \mbox{ and } k+l=3[\frac{n}{2}]+1\mbox{ or }=3[\frac{n}{2}],$ \\
and $x_{l+s,k+s}=x_{l,k}$ for all  $s\in \mathbb{Z}.$\\
{\bf For $n$ even}\\
$x_{l,k}=-(-1)^{k-l}{n\choose k-l }\mbox{ for }n>k>l\geq 0 \mbox{ and } k+l=\frac{n}{2}\mbox{ or }\frac{n}{2}-1,$\\
$x_{l,k}=(-1)^{l-k}{n\choose l-k }\mbox{ for }n> l>k\geq 0 \mbox{ and } k+l=3\frac{n}{2}\mbox{ or }3\frac{n}{2}-1,$\\
and $x_{l+s,k+s}=x_{l,k}$ for all  $s\in \mathbb{Z}.$\\

\noindent 
From the above one deduces for $0\leq k,\ell <n,\ k\neq \ell$ the
formulas:\\
 {\bf For $n$ odd} and $0\leq \ell <k$ one has $x_{\ell ,k}=-(-1)^{k-\ell}{n\choose k-\ell }$\\
 and $x_{\ell ,k}=-x_{k,\ell}$.\\ 
{\bf For $n$ even} and $0\leq \ell <k$ one has
$x_{\ell ,k}=-(-1)^{k-\ell}{n\choose k-\ell }$ if $k-\ell \leq \frac{n}{2} $ and 
$x_{\ell ,k}=(-1)^{k-\ell}{n\choose k-\ell }$ if $k-\ell > \frac{n}{2} $.\\
 {\bf  For $n$
even} and $0\leq k<\ell$ one has $x_{\ell ,k}=(-1)^{\ell -k}{n\choose \ell -k }$. \hfill $\square$ \end{theorem}

The second part of 3.1 is obtained from the first part by using
the equalities $x_{\ell,k}=x_{\ell',k'}$ if $\ell \equiv \ell',\ k\equiv k'$ modulo $n$ and the equalities
$x_{\ell ,k}=x_{\ell +s,k+s}$ for all $s\in \mathbb{Z}$.

\subsection{Generalised eigenvalues and formal monodromy}

The scalar operator $(z\frac{d}{dz})^n-z$ can be transformed into a matrix differential operator
$(z\frac{d}{dz}) + A$ where the entries of the matrix $A$ are in $\mathbb{C}(\{z^{-1}\})$. More precisely, the matrix $A$ has the form 
$\begin{pmatrix} 
0 & 1 & & & & \\
 & 0& 1& & & \\
 &  &  0&1 & & \\ 
& & &\ddots& \ddots& \\
& & & & 0&1\\
z& & & & & 0
\end{pmatrix}$. \\

In the case of $(z\frac{d}{dz})^n-z$, the differential module $\mathbb{C}((z^{-1/n}))\otimes M$ has a basis $b_0,\dots ,b_{n-1}$ such that $\delta b_j=-q_jb_j$ with $q_j=\zeta ^jz^{1/n}$ and $\zeta =e^{2\pi i/n}$. The $q_j$ are the generalized
eigenvalues and the matrix form of $\delta$, with respect to this basis, reads $z\frac{d}{dz}-diag(z^{1/n},\zeta z^{1/n},\dots ,\zeta ^{n-1}z^{1/n})$.\\

The symbolic solution space $V$ has the basis $\{e_j:=e^{\frac{1}{n}\zeta ^jz^{1/n}}b_j |\  j=0,\dots ,n-1\}$. The elements $b_j$ are unique up to multiplication by a constant. From the identities
 $\gamma V_q=V_{\gamma q}$ it follows that
these constants are choosen such that the formal monodromy $\gamma$ has the form $e_0\mapsto e_1\mapsto \cdots \mapsto e_{n-2}\mapsto e_{n-1}\mapsto (-1)^ne_0$. The sign $(-1)^n$ comes from the observation that $\gamma$ has determinant 1 on $V$.\\
 
In this case $mon_\infty$ can be identified with the topological monodromy $mon_0$ at $z=0$ (because $\mathbb{Z}$ is the fundamental group of $\mathbb{C}^*$). This is a unipotent matrix with characteristic polynomial $(\lambda -1)^n$. 
 
\subsection{The singular directions}

Put $(\zeta ^k-\zeta ^\ell)=| \zeta ^k-\zeta ^\ell |\cdot e^{2\pi i\phi (k,\ell)}$ with, say, $0\leq \phi(k,\ell) <1$. 
Now $d$ is a singular direction for $q_k-q_\ell$ if and only if  $\cos (2\pi \phi(k,\ell)+2\pi \frac{d}{n})=-1$. Thus $d=:d(k,l)=n(\frac{1}{2} -\phi(k,\ell))$ is modulo $n$ the only singular direction for the
pair $(q_k,q_l)$. 

Recall that the symbolic solution space $V$ has basis $e_0,\dots ,e_{n-1}$.
We denote by $E_{a,b}\in {\rm End}(V)$ the map given by $E_{a,b}e_b=e_a$ and $E_{a,b}e_c=0$
for $c\neq b$. One has $E_{a,b}E_{b,c}=E_{a,c}$. Moreover, the part of $St_{d(k,l)}$ corresponding to the
pair $(q_k,q_l)$ has the form $x_{l,k}E_{l,k}$ for a certain constant $x_{l,k}$. Then 
\[ St_d=1+\sum _{(k,l)\mbox{\small{ such that} }  d=d(k,l) \mbox{ \small{mod} } n}x_{l,k}E_{l,k}.\]
 {\it Our goal is to compute all constants $x_{l,k}$}.\\

\noindent {\it Computation of $d(k,l)$}.
One observes that for $\lambda \in (0,1)\subset \mathbb{R}$, the formula
$(e^{2\pi i \lambda}-1)=|(e^{2\pi i \lambda}-1)|e^{2\pi i \mu}$ holds with $\mu =\frac{1}{4}+\frac{\lambda }{2}$.
This implies:\\

\noindent For $n>k>l\geq 0$ one has $\phi(k,l)=\frac{1}{4}+\frac{k+l}{2n}$ and $d(k,l)=\frac{n}{4}-\frac{k+l}{2}$.\\
 For $n>l>k\geq 0$ one has $\phi(k,l)=\frac{3}{4}+\frac{k+l}{2n}$ and $d(k,l)=\frac{3n}{4}-\frac{k+l}{2}$.\\

\noindent  {\it  For $n$ odd and $n>k>l\geq 0$}, the possibilities for $d(k,l)\in [0,1)+\mathbb{Z}n$ are given by:
$k+l=[\frac{n}{2}],\ d(k,l)=\frac{1}{4}$ and $k+l=[\frac{n}{2}]-1,\ d(k,l)=\frac{3}{4}$.\\
 \noindent  {\it For $n$ odd and $n>l>k\geq 0$}, the possibilities for $d(k,l)\in [0,1)+\mathbb{Z}n$ are given by: $k+l=3[\frac{n}{2}]+1,\ d(k,l)=\frac{1}{4}$ and $k+l=3[\frac{n}{2}],\ d(k,l)=\frac{3}{4}$.\\

\noindent {\it  For $n$ even and $n>k>l\geq 0$}, the possibilities for $d(k,l)\in [0,1)+\mathbb{Z}n$ are given by: $k+l=\frac{n}{2},\ d(k,l)=0$ and $k+l=\frac{n}{2}-1,\ d(k,l)=\frac{1}{2}$.\\
 \noindent  {\it For $n$ even and $n>l>k\geq 0$}, the possibilities for $d(k,l)\in [0,1)+\mathbb{Z}n$ are given by: $k+l=3\frac{n}{2},\ d(k,l)=0$ and $k+l=3\frac{n}{2}-1,\ d(k,l)=\frac{1}{2}$.\\
 
\subsection{The equation for odd $n$}
  The monodromy identity for odd $n$ is:
  $mon_\infty \mbox{ is conjugated to } \gamma St_{\frac{3}{4}}St_{\frac{1}{4}}$. Therefore
  $P_n:=\det (-\lambda 1+ \gamma St_{\frac{3}{4}}St_{\frac{1}{4}})$ equals $-(\lambda -1)^n$. Further
  \[\gamma =E_{1,0}+E_{2,1}+\cdots +E_{n-1,n-2}+E_{0,n-1},\]
  \[St_{\frac{3}{4}}=1+\sum_{k+l=[\frac{n}{2}]-1, k>l}x_{l,k}E_{l,k} 
  +\sum_{k+l=3[\frac{n}{2}],l>k}x_{l,k}E_{l,k}, \]
  \[St_{\frac{1}{4}}=1+\sum_{k+l=[\frac{n}{2}], k>l}x_{l,k}E_{l,k} 
  +\sum_{k+l=3[\frac{n}{2}]+1,l>k}x_{l,k}E_{l,k}. \]
One observes that $P_n$ is the determinant of a sparse matrix
and guided by a few explicit examples, verified by a MAPLE,
 \[P_3=-\lambda ^3+x_{0,1}\lambda ^2+x_{2,1}\lambda +1,  \]
\[P_5=-\lambda ^5+x_{0,1}\lambda ^4+x_{0,2}\lambda ^3+x_{4,2}\lambda ^2+x_{4,3}\lambda +1,\]
\[P_7=-\lambda ^7+ x_{1,2}\lambda ^6+x_{0,2}\lambda ^5+x_{0,3}\lambda ^4+x_{6,3}\lambda ^3+x_{6,4}\lambda ^2+x_{5,4}\lambda +1,\]
\begin{small}\[P_9=-\lambda ^9+x_{1,2}\lambda^8+x_{1,3}\lambda^7+x_{0,3}\lambda^6+x_{0,4}\lambda^5+x_{8,4}\lambda^4+x_{8,5}\lambda^3+x_{7,5}\lambda^2+x_{7,6}\lambda+1.\]\end{small}
 one obtains {\bf
 the general formula for $P_n$ and odd $n$ }: $P_n=-\lambda ^n+1+$
\[\sum _{k>l,\ k+l=[\frac{n}{2}]\mbox{ or }=[\frac{n}{2}]-1}x_{l,k}\lambda ^{n-(k-l)}+
\sum _{l>k,\ k+l=3[\frac{n}{2}]+1\mbox{ or }=3[\frac{n}{2}]}x_{l,k}\lambda ^{l-k}.\]
From this and the equality $\gamma ^{-1}St_d\gamma =St_{d+1}$ one obtains
\[ x_{l,k}=-(-1)^{k-l}{n\choose k-l }\mbox{ for }k>l \mbox{ and } k+l=[\frac{n}{2}]\mbox{ or }=[\frac{n}{2}]-1, \]
\[x_{l,k}=(-1)^{l-k}{n\choose l-k }\mbox{ for } l>k \mbox{ and } k+l=3[\frac{n}{2}]+1\mbox{ or }=3[\frac{n}{2}], \]  
\ \ \ \  and  $x_{l+t,k+t}=x_{l,k}$ for all  $t\in \mathbb{Z}$.

The proof of the formula for $P_n$ consists simply of determining for each power of $\lambda$ the part of the sparse matrix which contributes to its coefficient in the determinant. The verification  is straightforward.
\subsection{The equation for even $n$}
According to the monodromy identity, $mon_\infty$ is conjugated to $\gamma St_{\frac{1}{2}}St_0$.
Thus $P_n:=\det (-\lambda 1+\gamma St_{\frac{1}{2}}St_0)$ equals $(\lambda -1)^n$.
\[\gamma =E_{1,0}+E_{2,1}+\dots +E_{n-1,n-2}-E_{0,n-1},\]
\[St_{\frac{1}{2}}=1+\sum _{k>l,\ k+l=\frac{n}{2}-1}x_{l,k}E_{l,k}+\sum_{l>k,\ k+l=3\frac{n}{2}-1}x_{l,k}E_{l,k}   ,\]
\[St_0=1+ \sum_{k>l,\ k+l=\frac{n}{2}} x_{l,k}E_{l,k}+\sum_{l>k,\ k+l=3\frac{n}{2}}x_{l,k}E_{l,k} .\]
Guided by a few examples, verified by a MAPLE computation
 \[P_2=\lambda ^2-x_{0,1}\lambda +1,\]
\[ P_4=\lambda ^4-x_{0,1}\lambda ^3-x_{0,2}\lambda ^2+x_{3,2}\lambda  +1,\]
\[P_6=\lambda ^6  -x_{1,2}\lambda^5-x_{0,2}\lambda ^4-x_{0,3}\lambda ^3+x_{5,3}\lambda ^2+x_{5,4}\lambda+1,\]
\begin{small}\[P_8=\lambda ^8 - x_{1,2}\lambda^7 -x_{1,3}\lambda ^6- x_{0,3} \lambda ^5- x_{0,4}\lambda ^4+ x_{7,4}\lambda ^3+ x_{7,5}\lambda ^2+ 
x_{6,5}\lambda +1.\]\end{small}
one deduces {\bf the general formula for $P_n$ and even $n$.} $P_n=\lambda ^n+1+$
\[-\sum _{k>l,\ k+l=\frac{n}{2}\mbox{ or }=\frac{n}{2}-1}x_{l,k}\lambda ^{n-(k-l)}+
\sum _{l>k,\ k+l=3\frac{n}{2}\mbox{ or }=3\frac{n}{2}-1}x_{l,k}\lambda ^{l-k}.\] 
This implies
\[ x_{l,k}=-(-1)^{k-l}{n\choose k-l }\mbox{ for }k>l \mbox{ and } k+l=\frac{n}{2}\mbox{ or }\frac{n}{2}-1, \]
\[x_{l,k}=(-1)^{l-k}{n\choose l-k }\mbox{ for } l>k \mbox{ and } k+l=3\frac{n}{2}\mbox{ or }3\frac{n}{2}-1, \]  
\ \ \ \ \ \ \ and $ x_{l+s,k+s}=x_{l,k}$ for all  $s\in \mathbb{Z}$. 

\begin{remark} Weighted projective spaces. \label{rm32} \\
{\rm Consider positive integers $w_0,\dots ,w_n$ with
$\gcd (w_0,\dots ,w_n)=1$.
For the weighted projective space $\mathbb{P}(w_0, .., w_n)$,
which is defined by $\mathbb{C}^{n+1}\setminus \{0\}/\mathbb{C}^*$,
where $t\cdot (z_0,\dots ,z_n)=(t^{w_0}z_0,\dots ,t^{w_n}z_n)$, we
adopt the  quantum differential operator, given in \cite{GuestSakai}, namely
\[ \prod _{i=1}^n(w_i\hbar \partial)(w_i\hbar \partial -\hbar)\cdots
(w_i\hbar -(w_i-1)\hbar)-q , \]
where $\partial=q\frac{d}{dq}$. After taking $\hbar =1$ and replacing 
$q$ by $z$ and $\partial$ by $\delta =z\frac{d}{dz}$ the operator reads
\[\prod _{j=0}^n\delta (\delta -\frac{1}{w_j})\dots (\delta -\frac{w_j-1}{w_j})\ -z.\]

We note that the above formula is attributed to Corti and  Golyshev and that in \cite{Iritani}  Dubrovin's conjecture is extended to orbifolds. In particular, there is a conjecture for weighted projective spaces. Unfortunately the latter is not explicit enough 
to allow us a comparison with the Stokes data.  Here we will show that
our computations of the classical Stokes matrices for ordinary projective
spaces extend to the case of weighted projective spaces. The preprint
\cite{T-U}, related to Proposition \label{p33},  appeared after this paper was finished. \\

\begin{proposition} \label{p33} {\bf The Stokes data $\{x_{\ell ,k}\}$ for} 
\[\prod _{j=0}^n\delta (\delta -\frac{1}{w_j})\dots (\delta -\frac{w_j-1}{w_j})\ -z. \mbox{ Put  } s=\sum w_j.\]
At $z=\infty$, the generalized eigenvalues are $\zeta ^j z^{1/s}$ with $j=0,\dots ,s-1$ where $\zeta=e^{2\pi i/s}$. 
Thus the above equation is formally equivalent to $\delta ^s-z$
 and the configuration of the Stokes matrices is the same as for the  ordinary projective space   $\mathbb{P}^{s-1}$.  The formal monodromy  differs by a minus-sign if $n$ is even. 

The topological monodromy at $z=0$ (or equivalently at $z=\infty$) has characteristic polynomial $\prod _{j=0}^n(\lambda ^{w_j}-1)$.

 The Stokes data $\{x_{\ell ,k}\}$ are determined by:\\
{\rm (a)}. The monodromy identity $\pm P_n= \prod _{j=0}^n(\lambda ^{w_j}-1)$. \\
{\rm (b)}. $x_{\ell ,k}=x_{\ell ',k'}$ if $\ell \equiv \ell ',\ k\equiv k' \mod s$.\\
{\rm (c)}. $x_{\ell +t,k+t}=x_{\ell ,k}$ for all $t\in \mathbb{Z}$.\\
In particular, the Stokes data consists of  computable integers.
 \end{proposition}
The proof is a straightforward computation. We note that it might be
difficult to give a closed formula (as in the $\mathbb{P}^{n-1}$ case)
for the $x_{\ell,k}$.\\

\noindent
{\it Example  $\mathbb{P}(1,2,4)$}.
The topological monodromy at $z=\infty$ is conjugated to 
$\gamma St_{3/4}St_{1/4}$. The characteristic polynomial of this $7\times 7$-matrix is\\
$-\lambda ^7+x_{1,2}\lambda^6+x_{0,2}\lambda ^5+x_{0,3}\lambda ^4+x_{6,3}\lambda ^3
+x_{6,4}\lambda ^2+x_{5,4}\lambda +1$, 
where these $x_{\ell,k}$ are the non trivial entries of $St_{3/4}$ and $St_{1/4}$. 

The topological monodromy at $z=0$ has characteristic polynomial\\
$-(\lambda-1)(\lambda ^2-1)(\lambda ^4-1)$ and thus we find\\ 
$x_{1,2}=1,\ x_{0,2}=1,\ x_{0,3}=-1 ,\ x_{6,3}=1 ,\ x_{6,4}=-1 ,\ x_{5,4}=-1$.  }\end{remark}

\section{The quantum differential equation \\
$\delta ^{N-1}-zk(k\delta +k-1)(k\delta +k-2)\cdots (k\delta +1)$} \label{s4}

According to \cite{Guest}, the Dubrovin--Givental connection 
for a non singular hypersurface of degree $k\leq N-1$
in $\mathbb{P}^{N-1}$ is given by this formula.
  We prefer to write  this operator differently (with
$m=k$ and $n=N-k$) \begin{small}
\[\delta ^{n+m-1}-m^mz(\delta +\frac{m-1}{m})(\delta +\frac{m-2}{m})\cdots (\delta +\frac{1}{m})
\mbox{ with } \delta=z\frac{d}{dz},\ n>1, m>1 .\] \end{small} 
For $m=1$ this reduces this operator to the one studied in \S 3. At the end of this section we will comment on the case $n=1$. 

\begin{theorem} \label{teorema41} The Stokes data for the above equation is:\\ $\{x_{\ell,k}|\ 0\leq k, \ell \leq n-1,\  k\neq \ell\}$ and $\{z_j|\ 1\leq j \leq m-1\},\ \{y_j|\ 1\leq j \leq m-1\}$. The $y_j$ and $z_j$ depend on the choice of
a basis. However, the products $y_jz_j,\ j=1,\dots , m-1$ are computable elements of $\mathbb{Q}(\zeta )$, where $\zeta=e^{2\pi i/m}$ and independent of this choice. \\ 
 $x_{\ell+s,k+s}=x_{\ell,k}$ holds for $s\in \mathbb{Z}$ and
$x_{\ell ,k}=x_{\ell ',k'}$ if $\ell \equiv \ell ', \ k\equiv k' \mod n$.\\
{\bf For $n>1$ odd}\\
$ x_{\ell ,k}=(-1)^{k-\ell +1}{n+m\choose k-\ell }$ for $n>k>\ell \geq 0 ,\ k+\ell=[\frac{n}{2}]$ or $=[\frac{n}{2}]-1$,\\
$ x_{\ell , k}=(-1)^{\ell -k}{ n+m\choose  \ell -k}$ for $n>\ell >k\geq 0,\ k+\ell =3[\frac{n}{2}]+1$ or $=3[\frac{n}{2}].$ \\ 
{\bf For $n$ even}\\
$x_{\ell ,k}=(-1)^{k-\ell +1}{n+m\choose k-\ell }$ for $n>k>\ell \geq 0, k+\ell =\frac{n}{2}$ or $=\frac{n}{2}-1$,  \\
$ x_{\ell ,k}=(-1)^{n+m+\ell -k +1}{n+m\choose \ell -k }$ for $n>\ell> k\geq 0, k+\ell =3\frac{n}{2}$ or $=3\frac{n}{2}-1$.  
\end{theorem}

\subsection{The differential equation $\delta ^4-27z\delta ^2-27z\delta -6z$}
We start by investigating the case $n=2, m=3$ of Theorem 4.1,
which is the quantum differential equation of a hypersurface of degree 3 in $\mathbb{P}^4$ (see \cite{Guest}, p 42, Example 3.6).
 A matrix form for this equation is
\[z\frac{d}{dz}+\left(\begin{array}{cccc} 0&-1 &0 &0 \\ 0 &0 &-1 & 0 \\ 0 &0 & 0&-1 \\
-6z &-27z & -27z &0 \end{array}\right).\]
We proceed as in \S 3. The (generalised) eigenvalues at $z=\infty$ are 
$q_1=\sqrt{27}z^{1/2}, \ q_2=-\sqrt{27}z^{1/2},\ 0$. The symbolic solution space $V$ at $z=\infty$
has the form $V=V_{q_1}\oplus V_{q_2}\oplus V_0$ with $V_{q_1}=\mathbb{C}e_1$,
$V_{q_2}=\mathbb{C}e_2$ and $V_0=\mathbb{C}e_3\oplus \mathbb{C}e_4$. The basis 
$e_1,\dots ,e_4$ is chosen such that the formal monodromy has the form
\[\gamma =\left( \begin{array}{cccc}0 &-1 & 0& 0\\ 1& 0 & 0 & 0\\ 0& 0& \zeta & 0\\
0 &0 &0  & \zeta ^2 \end{array}\right), \mbox{ where } \zeta =e^{2\pi i/3}. \]
We note that this basis is unique up to a transformation of the type\\  $e_1\mapsto \lambda _1e_1,\
e_2\mapsto \lambda _1e_2,\  e_3\mapsto \lambda _2e_3,\ e_4\mapsto \lambda _3e_4$ with all
$\lambda _j\in \mathbb{C}^*$.

The singular directions are $0+2\mathbb{Z}$ for the differences $q_2-q_1,\ q_2-0,\ 0-q_1$ and are
$1+2\mathbb{Z}$ for the differences $q_1-q_2,\ q_1-0,\ 0-q_2$. The Stokes matrix $St_0$ has the form
\begin{small}
\[ St_0=\left(\begin{array}{cccc} 1& 0&x_4 & x_5 \\ x_1 &1 & 0& 0\\  0& x_2 & 1 & 0 \\ 0 & x_3 & 0 & 1\end{array}\right),\  \gamma St_0=\left(\begin{array}{cccc} -x_1& -1& 0 & 0 
\\ 1 &0 & x_4& x_5\\  0& \zeta x_2 & \zeta & 0 \\ 0 & \zeta ^2x_3 & 0 & \zeta ^2\end{array}\right) \] \end{small}
and $ St_1=\gamma ^{-1} St_0\gamma $.
According to the monodromy identity, $\gamma St_0$ is equivalent to the topological monodromy
at $z=0$. The latter is seen to have the single eigenvalue 1 (and only one Jordan block).
Thus the characteristic polynomial of $\gamma St_0$ is $(\lambda -1)^4$. This yields the data for the
entries of the Stokes matrices $ x_1=-5,\ x_2x_4=-9\zeta +18,\ x_3x_5=9\zeta  +27$. It seems that we have
not enough information to obtain values for all $x_j$. This is due however to the non uniqueness of the basis vectors $e_3,e_4$. As an example we can see that for a suitable choice of $e_3,e_4$ we will have, say, $x_4=1$ and $x_5=1$
and further $x_1=-5,\  x_2=-9\zeta +18,\ x_3 =9\zeta +27$. 

\subsection{The general case}
The above operator is transformed in the usual way into a first order matrix differential operator.
The formal data for the symbolic solution space $V$ at $z=\infty$ are: the (generalised)
eigenvalues are $0$ and the  $q_j=\sqrt[n]{m^m}\zeta _n^jz^{1/n}$ for $j=0,1,\dots ,n-1$ with 
$\zeta _n=e^{2\pi i/n}$. This solution space $V$ has the decomposition
$V=V_{q_0}\oplus V_{q_1}\oplus \cdots \oplus V_{q_{n-1}}\oplus V_0$ with 
$V_{q_j}=\mathbb{C}e_j$ for $j=0,\dots ,n-1$ and $V_0=\mathbb{C}f_1\oplus \cdots \oplus \mathbb{C}f_{m-1}$.  The basis vectors are chosen such that the formal monodromy
$\gamma$ acts as $e_0\mapsto e_1\mapsto \cdots \mapsto e_{n-1}\mapsto  (-1)^{n-1}(-1)^{m-1}e_0$  
and $\gamma f_j=\zeta _m^jf_j$ for $j=1,\dots ,m-1$ and $\zeta _m=e^{2\pi i/m}$.\\
We note that the basis  $f_1,\dots ,f_{m-1}$ of $V_0$ is unique up to multiplication by scalars.
The computation of the `monodromy identity' is done separately for $n$ even and $n$ odd.\\

\noindent {\bf even $n$}.\\
The singular directions $d$ for $q_k-q_\ell$ lying in $[0,1)+\mathbb{Z}n$ are the same as in \S 3,  namely:\\  
For $n>k>\ell \geq 0$: $d=0$ and $k+\ell =\frac{n}{2}$; $d=\frac{1}{2}$ and $k+\ell =\frac{n}{2}-1$.\\
For $n>\ell >k\geq 0$: $d=0$ and $k+\ell=\frac{3n}{2}$; $d=\frac{1}{2}$
 and $k+\ell =\frac{3n}{2}-1$.\\

\noindent For $q_k-0$, the only singular direction in $[0,1)+\mathbb{Z}n$ is $d=0$ with $k=\frac{n}{2}$.\\
 For $0-q_k$, the only singular direction in $[0,1)+\mathbb{Z}n$ is $d=0$ with $k=0$.\\

\noindent {\it  Description of  $St_0$}. 
For elements in ${\rm End}(\mathbb{C}e_0+\cdots +\mathbb{C}e_{n-1})$ we use the notation of \S 3. Then $St_0$ is the identity plus a number of maps, namely 
$\sum _{k>\ell, \ k+\ell =\frac{n}{2}} x_{\ell ,k}E_{\ell, k}$ and 
$\sum _{\ell >k,\ k+\ell =\frac{3n}{2}} x_{\ell ,k}E_{\ell, k}$ and a map 
$e_{\frac{n}{2}}\mapsto y_1f_1+\cdots +y_{m-1}f_{m-1}$ (the other base vectors are mapped to 0)
and for $j=1,\dots ,m-1$ a map $f_j\mapsto z_je_0$ (the other base vectors are mapped to 0).\\

\noindent {\it  Description of  $St_{\frac{1}{2}}$}. 
This Stokes matrix is the identity plus certain maps, namely $\sum _{k>\ell ,\ k+\ell =\frac{n}{2}-1}
x_{\ell ,k}E_{\ell, k}$ and $\sum _{\ell >k,\ k+\ell = \frac{3n}{2}-1}x_{\ell ,k}E_{\ell,k}$.\\

The matrix $\gamma St_{\frac{1}{2}}St_0$ and its characteristic polynomial  $P$ can be computed.
The monodromy identity $P=(\lambda -1)^{n+m-1}$ leads to the statement that
$x_{\ell,k}$ have the form $\pm {n+m\choose * }$ and that the  $y_jz_j$ are  elements of
  $\mathbb{Q}[\zeta _m]$. As in \S 4.1, i.e., the case $n=2,m=3$, one cannot compute
  the $y_j$ and $z_j$ separately since this involves a definite choice of the basis $f_1,\dots ,
  f_{m-1}$.\\
  
  \noindent {\it Example}. The case $n=4$, $m=3$ and $\zeta :=e^{2\pi i /3}$. \begin{small}
\[ \gamma =\left(\begin{array}{cccccc} 0& 0&0 &-1 &0 &0 \\ 1 &0 &0 & 0& 0&0\\ 
0&1 &0 & 0& 0&0 \\  0& 0& 1 &0 &0 & 0 \\ 0& 0 & 0& 0& \zeta & 0 \\ 0 & 0 & 0 & 0 & 0 & \zeta ^2     \end{array}\right),\
St_{\frac{1}{2}}= \left(\begin{array}{cccccc} 1& x_{0,1} &0 &0 & 0& 0 \\ 0 &1 &0 &0 &0 & 0\\ 
0& 0 & 1 & 0& 0& 0\\ 0 & 0 & x_{3,2} & 1 &0 &0 \\  0& 0& 0& 0& 1 &0\\ 0 &0 &0 &0 &0 & 1    
 \end{array}\right),\]
\[St_0= \left(\begin{array}{cccccc} 1& 0& x_{0,2}& 0& z_1 &z_2 \\
0 &1 & 0& 0& 0&0 \\ 0&0 &1 & 0& 0&0 \\ 0&0 &0 & 1& 0&0 \\ 0 &0 &y_1 & 0 &1 & 0
\\ 0& 0& y_2& 0& 0& 1    \end{array}\right).\]\end{small}
Since the characteristic polynomial of $\gamma St_{\frac{1}{2}}St_0$ is $(\lambda -1)^6$ one finds\\
 $x_{0,1}=7,\ x_{0,2}=-21,\ x_{3,2}=-7$, $y_1z_1= 9(2\zeta ^2+1)  $, $y_2z_2= -9(2\zeta ^2+1) $.\\

\noindent Let $P$ denote again the characteristic polynomial  of $\gamma St _{\frac{1}{2}}St_0$ for $n$ even and $m>1$.
One observes that $(\lambda -1)^{n+m}=(\lambda -1)P$ is the sum of 
$(\lambda ^m-1)Q$ with 
\[Q=\lambda ^n-\sum _{k>l,\ k+l=\frac{n}{2}\mbox{ or }=\frac{n}{2}-1}x_{l,k}\lambda ^{n-(k-l)}+
\sum _{l>k,\ k+l=3\frac{n}{2}\mbox{ or }=3\frac{n}{2}-1}x_{l,k}\lambda ^{l-k} +1 \]
and terms $a\lambda ^j$ ($a\in \mathbb{C}$) with $3\frac{n}{2}<j< m+\frac{n}{2}$. This leads to the formulas
\[ x_{\ell ,k}=(-1)^{k-\ell +1}{n+m\choose k-\ell }\mbox{ for } k>\ell, k+\ell =\frac{n}{2}\mbox{ or }=\frac{n}{2}-1,  \]
\[ x_{\ell ,k}=(-1)^{n+m+\ell -k +1}{n+m\choose \ell -k }\mbox{ for } \ell> k, k+\ell =3\frac{n}{2}\mbox{ or }=3\frac{n}{2}-1.  \]
The elements $y_jz_j$ are (in general complicated) expressions in $\mathbb{Q}(\zeta )$.\\

\noindent {\bf odd $n>1$}.\\

{\it The singular directions $d$ in $[0,1)+\mathbb{Z}n$ are}:\\
For $q_k-q_\ell$:\\
$n>k>\ell \geq 0$, $d=\frac{1}{4}$ with $k+\ell =[\frac{n}{2}]$; $d=\frac{3}{4}$ with 
$k+\ell =[\frac{n}{2}]-1$\\
$n>\ell >k\geq 0$, $d=\frac{1}{4}$ with $k+\ell =3[\frac{n}{2}]+1$; $d=\frac{3}{4}$ with
 $k+\ell=3[\frac{n}{2}]$.\\
For $q_k-0$: $d=\frac{1}{2}$ and $k=[\frac{n}{2}]$.
For $0-q_k$: $d=0$ and $k=0$.\\

\noindent {\it Example. $n=3,m=3$, $\zeta =e^{2\pi i/3}$}.\begin{small}
\[ \gamma:=\left(\begin{array}{ccccc}  0&0 &  1& 0&0 \\ 1 &0 &0 & 0&0\\ 0 &1 &0 &0 &0\\
  0&0 &0 &\zeta & 0\\ 0 &0 &0 &0 & \zeta ^2    \end{array}\right),\
 St_{\frac{3}{4}}=\left(\begin{array}{ccccc}  1&0 &0 & 0&0 \\ 0 &1 &0 &0 &0 \\ 0 &x_{2,1} &1 &0 &0\\  0&0 &0 &1 &0\\ 0 &0 &0 &0 & 1   \end{array}\right),\]
\[ St_{\frac{1}{2}}=\left(\begin{array}{ccccc} 1 &0 & 0& 0& 0\\ 0 &1 &0 &0 &0\\ 0 &0 &1 &0 &0\\ 
0 &z_1 &0 &1 &0\\  0&z_2 &0 &0 & 1   \end{array}\right), \
 St_{\frac{1}{4}}=\left(\begin{array}{ccccc} 1 &x_{0,1} &0 &0 &0 \\ 0 &1 &0 &0 &0\\
  0 &0 &1 &0 &0\\ 0 &0 &0 &1 &0\\  0& 0& 0& 0& 1   \end{array}\right),\]
\[ St_0=\left(\begin{array}{ccccc}  1&0 &0 &y_1 &y_2 \\ 0 &1 &0 &0 &0\\ 0 &0 &1 &0 &0\\ 
0 &0 &0 &1 &0\\  0& 0& 0& 0& 1   \end{array}\right).\] \end{small}

The observation that the characteristic polynomial of $\gamma St_{\frac{3}{4}}St_{\frac{1}{2}}
St_{\frac{1}{4}}St_0$ is $(\lambda -1)^5$ yields $x_{0,1}=6, \ x_{2,1}=-6$ and 
$y_1z_1=-9(\zeta ^2+1),\ y_2z_2=9\zeta ^2$. \\

As in the case of even $n$ one obtains  {\it  for general odd $n>1$ and $m>1$} explicit formulas for
the entries  $x_{\ell,k}$ (same notation as in the even case) of the Stokes matrices, namely
\[ x_{\ell ,k}=(-1)^{k-\ell +1}{n+m\choose k-\ell }\mbox{ for } k>\ell ,\ k+\ell=[\frac{n}{2}]\mbox{ or }=[\frac{n}{2}]-1,\]
\[ x_{\ell , k}=(-1)^{\ell -k}{ n+m\choose  \ell -k}\mbox{ for } \ell >k,\ k+\ell =3[\frac{n}{2}]+1\mbox{ or }=3[\frac{n}{2}]  .\] 
The elements $y_jz_j$ are (in general complicated) expressions in $\mathbb{Q}(\zeta )$.\\

\noindent {\it Comments on the case $n=1$}.\\
The equation reads
$\delta ^{m}-m^mz(\delta +\frac{m-1}{m})(\delta +\frac{m-2}{m})\cdots (\delta +\frac{1}{m})$. The (generalized) eigenvalues at $z=\infty$ are $z$ and $0$. 
{\it This equation is not really a quantum differential equation and moreover there is no
ramification  at $z=\infty$}!

The symbolic solution space $V$ is given
a basis $e_0,f_1,\dots ,f_{m-1}$ such that $V_{z}=\mathbb{C}e_0$, $V_0$ has basis
$f_1,\dots ,f_{m-1}$ and the formal monodromy $\gamma$ has the form $\gamma (e_0)=e_0$
and $\gamma (f_j)=\zeta ^jf_j$ for all $j$ and $\zeta =e^{2\pi i/m}$. The above basis is unique up to
multiplication by scalars. The singular directions are $d=\frac{1}{2}$ and $d=0$ and the corresponding
Stokes matrices involve (using the earlier notation) only $\{y_1,\dots ,y_{m-1}\}$ and $\{z_1,\dots ,z_{m-1}\}$. These elements are not unique, however the products $y_jz_j$ are independent of the choice of $e_0,f_1,\dots ,f_{m-1}$ and are computable elements of $\mathbb{Q}(\zeta)$. \\
{\it Example}:  for $m=3$ one finds  $y_1z_1=3+3\zeta$, $y_2z_2=-3\zeta$. 
This example seems unrelated to the quantum cohomology of a cubic surface, studied by K. Ueda in \cite{Ueda2}.\\

\textbf{Acknowledgements} The first author would like to thank  Hiroshi Iritani for pointing out the reference \cite{Tan} and for many interesting conversations about the quantum cohomology of Fano varieties. JACM was supported by a Japanese Government (Monbukagakusho:MEXT) Scholarship. We thank the referees for their
helpful comments.

\begin{small}

\end{small}

\end{document}